\def\beqnas{\begin{eqnarray*}}
\def\eeqnas{\end{eqnarray*}}
\def\beqna{\begin{eqnarray}}
\def\eeqna{\end{eqnarray}}

\def\intl{\int\limits}
\def\suml{\sum\limits}
\def\txt#1{\quad\mbox{#1}\quad}

\def\myJ{J_{W}}
\def\myJhat{\widehat{J}_{W}}
\def\Jfnl{\mathbbm{J}}
\def\myG{\mathbbm{G}}
\def\myF{\mathbbm{F}}

\def\E{\mathbbm{E}}
\def\P{\mathbbm{P}}
\def\R{\mathbbm{R}}
\def\1{\mathbbm{1}}
\def\phi{\varphi}

\documentclass[a4paper,12pt]{article}
\usepackage{amssymb, amsthm, latexsym, bbm, graphicx, natbib}
\usepackage[latin1]{inputenc}
\usepackage{anysize}
\marginsize{2cm}{2cm}{2cm}{2cm}

\newcommand{\dee}[1]{\,{\rm d}{#1}}	

\makeatletter
\renewcommand\section{\@startsection{section}{1}{\z@}%
                                   {-3.5ex \@plus -1ex \@minus -.2ex}%
                                   {2.3ex \@plus.2ex}%
                                   {\normalfont\bfseries}}
\renewcommand\subsection{\@startsection{subsection}{2}{\z@}%
                                     {-3.25ex\@plus -1ex \@minus -.2ex}%
                                     {1.5ex \@plus .2ex}%
                                     {\normalfont\itshape}}
\renewcommand\th@plain{%
  \thm@headfont {\scshape}
  \itshape}
\renewcommand\th@definition{%
  \thm@headfont {\scshape}
  \normalfont}

\makeatother

\begin{document}

\thispagestyle{empty}

\def\propositionname{Proposition}
\theoremstyle{plain}
\newtheorem{proposition}{\propositionname}
\theoremstyle{definition}
\newtheorem{definition}{Definition}
\def\proofname{Proof}
\def\pf{\ifvmode\else\newline\fi\noindent\textsc{\proofname:\ }}
\def\qed{\mbox{ $\Box$}}

\author{A.\ Baddeley\footnote{{\tt adrian@maths.uwa.edu.au}}\,,
        M.\ Kerscher$^{\ddag}\!$,  K.\ Schladitz$^{+}$, B.\ T.\ Scott\\
        {\small\sl Department of Mathematics and Statistics, 
                   University of Western Australia,}\\[-3pt]
        {\small\sl Nedlands WA 6907, Australia}\\[3pt]
        {\small\sl $^{+}$ Institut f\"ur Techno- und
	Wirtschaftsmathematik,}\\[-3pt]
        {\small\sl Erwin-Schr\"odinger-Stra\ss e, 
		   D-67663 Kaiserslautern, Germany}\\[3pt]
	{\small\sl $^{\ddag}$ Sektion Physik, 
	Universit\"at M\"unchen,}\\[-3pt] 
        {\small\sl Theresienstr. 37/III, 80333 M\"unchen, Germany}\\[0.3cm]}
\date{}
\title{Estimating the $J$ function without \\ edge 
        correction\thanks{Adrian Baddeley and Katja Schladitz
        were supported by a grant from the Australian Research Council. Martin 
        Kerscher was supported by the ``Sonderforschungsbereich SFB 375 f\"ur 
        Astroteilchenphysik der Deutschen Forschungsgemeinschaft''. Bryan Scott 
	was supported by an Australian Postgraduate Award and
        CSIRO Postgraduate Studentship.}}
\maketitle

\begin{abstract}
\vspace{-30pt}
  \noindent
  The interaction between points in a spatial point process 
  can be measured by its empty space function $F$, 
  its nearest-neighbour distance distribution function $G$,
  and by combinations such as the J-function $J = (1-G)/(1-F)$.
  The estimation of these functions is hampered by edge effects:
  the uncorrected, empirical distributions of distances observed in
  a bounded sampling window $W$ give severely biased estimates
  of $F$ and $G$. However, in this paper we show that 
  the corresponding {\em uncorrected\/} estimator of the function 
  $J=(1-G)/(1-F)$
  is approximately unbiased for the Poisson case, and is useful
  as a summary statistic. 
  Specifically, consider the estimate
  $\myJhat$ of $J$ computed from uncorrected estimates of $F$ and $G$.
  The function $\myJ(r)$, estimated by $\myJhat$, 
  possesses similar properties to the $J$ function, for example
  $\myJ(r)$ is identically $1$ for Poisson processes.  This enables
  direct interpretation of uncorrected estimates of $J$, something
  not possible with uncorrected estimates of either $F$, $G$ or
  $K$.  We propose a Monte Carlo test for complete spatial randomness
  based on testing whether $\myJ(r)\equiv 1$. Computer simulations suggest this 
  test is at least as powerful as tests based on edge corrected estimators
  of $J$.\\

\vspace{30pt}
\end{abstract}

\newpage
\section{Introduction}

A spatial point pattern is often studied by estimating
point process characteristics such as the empty space function $F$,
the nearest-neighbour distance distribution function $G$,
and the $K$-function. Here we consider 
\begin{equation}
  \label{defj}
  J(r)=\frac{1-G(r)}{1-F(r)}
\end{equation}
as advocated by \cite{badd:lies96}. The $J$ function is identically
equal to 1 for a Poisson process, and values of $J(r)$ less than or greater
than 1 are suggestive of clustering or regularity, respectively. (Note 
that it is of course possible to find non-Poisson point 
processes for which $J(r)=1$, as \cite{bedf:berg97} have shown). 

In practice, observation of the point process is usually restricted to some 
bounded window $W$. As a consequence, estimation of the summary functions,
which is based on the measurement of various distances of the point
process, is hampered by the ``edge effects'' (bias and censoring)
introduced by restricting observation of these distances 
to $W$. 
In order to counter edge effects it is necessary to apply
some form of edge correction to the empirical estimates of the summary
functions. For further details see \cite{badd98, skm89, cres91, ripl88}.


Unbiasedness is highly desirable when a summary function estimate is 
to be compared directly to the corresponding theoretical value for a
point process model. However, as Diggle argues in the discussion
of \cite{ripl77} and in \cite{digg83}, 
unbiasedness is not essential when using a summary function estimator as the 
test statistic in a hypothesis test,  since the bias will be accounted for
in the null distribution of the test statistic.

This paper studies the {\em uncorrected\/} estimator of $J$ 
obtained by ignoring edge effects and computing
$\widehat J = (1 - \widehat G)/(1 - \widehat F)$ from
the uncorrected, empirical distributions 
$\widehat G$ and $\widehat F$ of distances observed in a compact window.
It was prompted by the accidental discovery that this uncorrected estimator
remarkably still yields values $\widehat J(r)$ approximately equal to 1
for the Poisson process. 
An intuitive explanation is that the relative bias due to edge effects is
roughly equal for the estimates of $1-G$ and $1-F$, so that these biases
approximately cancel in the ratio estimator of $J$. 
It follows that the uncorrected estimate of $J$ could be used
for the direct visual assessment of deviations from the Poisson process,
something not possible with the uncorrected estimates of $F$, $G$ or $K$.
Our aim is to formalise this uncorrected estimator of $J$,
and to investigate its use as a summary statistic and as a test statistic
in point pattern analysis.

The paper is organised as follows. Section 2 outlines and generalises some 
fundamental point process ideas. In Section 3 we define the $\myJ$ function 
estimated by the uncorrected procedure described above and derives some of 
its properties. In Section 4 we verify that the natural estimator of $\myJ$ is 
the uncorrected estimate $\myJhat$. Finally, Section 5 presents
the results of a computational experiment to compare the
power of Monte Carlo Tests constructed from estimates of $J$ and
$\myJ$ as well as estimation results from the simulations of various point
process models in square and rectangular windows in $\R^2$ and in a
cubic window in $\R^3$.

\section{Background}

Let $X$ be a stationary point process in $\R^d$ with intensity
$\lambda$. \cite[For details of the theory of point processes 
see][]{dale:vere88,cres91,skm89}.
The empty space function $F(r)$ is the probability of finding a point
of the process within a radius $r$ of an arbitrary fixed point:
\begin{equation}
  \label{deff}
  F(r)=\P(X\cap B(0,r)\neq\emptyset),
\end{equation}
where $B(x, r)$ denotes the ball of radius $r$ centred at $x$.
The nearest neighbour distance distribution function $G(r)$ is the
probability of finding another point of the process in the ball of
radius $r$ centred at a ``typical'' point of the process:
\begin{equation}
  \label{defg}
  G(r)=\P^{0!}(X\cap B(0,r)\neq\emptyset),
\end{equation}
where $\P^{0!}$ denotes the reduced Palm distribution at the origin
$0$. Roughly speaking $\P^{0!}$ is the distribution of rest of the process
$X\setminus\{0\}$ given there is a point of the process at the
origin \cite[see][]{dale:vere88}.

Let $W$ be a compact observation window in $\R^d$ with nonempty interior.
The construction of estimators for $F$ is based on the stationarity 
of $X$ yielding
\begin{equation}
  \label{estf}
  F(r)=\frac1{|W|}\intl_W\P(X\cap B(x, r)\neq\emptyset)\,dx,
\end{equation}
where $|W|$ denotes the volume of $W$. For $G$ the Campbell-Mecke formula 
\cite[(4.4.3)]{skm89} gives
\begin{equation}
  \label{estg}
  G(r)=\frac1{\lambda|W|}\E\sum_{x\in X\cap W}
       \1\{X\cap B(x, r)\setminus\{x\}\neq\emptyset\}.
\end{equation}
In both cases we need information about $X\cap B(x, r)$ for $x\in W$, whereas
we only observe $X\cap B(x, r)\cap W$. Usually this edge effect problem is 
countered by restricting the integration in (\ref{estf}) and summation in
(\ref{estg}) to those points $x$ for which $B(x, r)\subseteq W$ (the ``border 
method'') or by weighting the contributions to the integral and sum so as to
correct for the bias \cite[see for example][]{badd98}.

The uncorrected, empirical distributions of distances observed in the 
window $W$ correspond to simply replacing $X$ by $X\cap W$ in (\ref{estf}) and
(\ref{estg}). In order to investigate the effect of this, we extend $F$ and $G$ to 
functionals as follows.

\begin{definition}
  \label{newfg}
  For every compact set $K\subset\R^d$ containing the origin define
\begin{eqnarray*}
  \myF(K) &:=& \P(X\cap K\neq\emptyset) \\
  \myG(K) &:=& \P^{0!}(X\cap K\neq\emptyset) \\
  \Jfnl(K) &:=& \frac{1-\myG(K)}{1-\myF(K)}.
\end{eqnarray*}
\end{definition}
\noindent
(Note that the empty space functional bears some relation to the contact
distribution function \citep[p.\ 105]{skm89} in that $H_{B}(r)=\myF(rB)$).
>From these we are able to define the window based $J$ function.

\section{The $\myJ$ function}

\begin{definition}
 For every compact set $W\subset\R^d$ with nonempty interior let
 \begin{equation}
   \label{myJ}
   \myJ(r):=\frac{\intl_W [1-\myG(B(0,r)\cap W_{-x})]\,dx}
                    {\intl_W [1-\myF(B(0,r)\cap W_{-x})]\,dx}
 \end{equation}
 be the \emph{window based $J$ function}, where $W_{-x}=\{y-x:y\in W\}$ is the
 translate of $W$ by $-x\in\R^d$.
 \label{myJdefn}
\end{definition}

If $X$ is a stationary Poisson process, then by
Slivnyak's Theorem \cite[(4.4.7)]{skm89} $\P\equiv\P^{0!}$. Thus 
$\myF\equiv\myG$ and we arrive at the following proposition. 
\begin{proposition}
  \label{newjpois}
  Let $X$ be a stationary Poisson process. Then 
  \[\myJ(r)\equiv 1\txt{for all}W\txt{and}r\geq0.\]
\end{proposition}

Explicit evaluation of $\myJ$ for other point process models seems difficult.
However, we can show that $\myJ$ behaves similarly
to the $J$ function for ordered and clustered processes,
suggesting it can also be interpreted in the same way as the $J$ function.
For some processes it is also possible to demonstrate that the $\myJ$ function
exhibits less deviation from the Poisson hypothesis than the equivalent
$J$ function.

\begin{proposition}
  \label{bound}
  Suppose $X$ is a process which is ``ordered'' in the sense that its
  $J$ functional is non-decreasing, that is $K_1\subseteq K_2$ implies
   $\Jfnl(K_1)\leq \Jfnl(K_2)$. Then
   \[ 
	1 \le \myJ(r) \le  J(r)
	\quad\mbox{for all } r. 
   \]
   Similarly if a process is ``clustered'',
   $K_1\subseteq K_2$ implies $\Jfnl(K_1) \ge \Jfnl(K_2)$,
   then $1 \ge \myJ(r) \ge J(r)$ for all $r$.
\end{proposition}
\pf
Observe that (\ref{myJ}) can be rewritten
\begin{equation}
  \myJ(r)
         =\int_{W}\Jfnl(B(0,r)\cap W_{-x})h_{W,r}(x)dx \label{weight}
\end{equation}
where 
\[
	h_{W,r}(x) = \frac{
			\left(
				1-\myF(B(0,r)\cap W_{-x})
			\right)
		  }{
			\int_W
			\left(
				1-\myF(B(0,r)\cap W_{-y})
			\right)
			\dee y
		  }
\]
satisfies $h_{W,r}(x) \geq 0$ for all $x\in \R^{d}$ 
and $\int_{W}h_{W,r}(x)dx=1$.
Hence
\[ \min_{W}\Jfnl(B(0,r)\cap W_{-x})\leq\myJ(r)\leq 
    \max_{W}\Jfnl(B(0,r)\cap W_{-x}) \]
and since $\Jfnl$ is nondecreasing
\[ 1= \Jfnl(\emptyset) \leq \myJ(r) \leq \Jfnl(B(0,r)) = J(r). \ \mbox{\qed}\]

The latter result can be strengthened to strict inequality
for specific examples. 
		
\begin{proposition}
  \label{PNS}
  Let $X$ be a Neyman-Scott cluster process with mean number of points
  per cluster greater than $1$. Assuming the support of the
  distribution of the cluster points contains a neighbourhood of the
  origin, then
  \[J(r)<\myJ(r)<1\txt{for all}W\txt{and}r>0.\]
\end{proposition}

Examples of processes satisfying the conditions of Proposition \ref{PNS}
are Mat\'ern's cluster process and the modified 
Thomas process described, for example, in \cite{skm89}. 

\pf
  The Palm distribution of a Neyman-Scott process is the convolution
  of the original distribution $P$ of the process and the Palm
  distribution $c_0$ of the representative cluster $N$
  \cite[(5.3.2)]{skm89}. Thus, for every compact set $K$, we have
  \beqnas
  1-\myG(K)&=&\intl\intl\1\{(\phi\cup\psi)\cap K\setminus\{0\}
              =\emptyset\}c_0(d\psi)P(d\phi)\\
           &=&\intl\intl\1\{\phi\cap K\setminus\{0\}=\emptyset\}
                        \1\{\psi\cap K\setminus\{0\}=\emptyset\}
                        c_0(d\psi)P(d\phi)\\
           &=&\P(X\cap K\setminus\{0\}=\emptyset)
              c_0(N\cap K\setminus\{0\}=\emptyset)\\
           &=&(1-\myF(K))c_0(N\cap K\setminus\{0\}=\emptyset).
  \eeqnas
Hence
\[\Jfnl(K) = c_0(N\cap K\setminus\{0\}=\emptyset). \]
  Now the assumption on the cluster distribution ensures 
  $c_0(N\cap K\setminus\{0\}=\emptyset)<1$ for all $K$ containing a
  neighbourhood of the origin. The conclusion then follows since $B(0,r)\cap W_{-x}$
  contains a neighbourhood of $0$ whenever $x$ is an interior point of $W$.
  \ \qed

\begin{proposition}
  Let $X$ be a hard-core process with hard-core radius $R$. Then
  \[1<\myJ(r)<J(r)\txt{for all}W\txt{and}0<r<R,\]
  and $J_{W}(r)$ is non-decreasing in $r$ for all $r<R$.
\end{proposition}
\pf
  Trivially we have $\myG(K)=0$ for all $K\subset B(0,R)$, while for
any point process the empty space functional $\myF$ is 
nondecreasing. Therefore the $J$ functional becomes 
\[\Jfnl(K) = \frac{1}{1-\myF(K)}\txt{for all}K\subset B(0,R),\]
which is also non-decreasing and the result follows by proposition \ref{bound}.
\ \qed

\section{Estimation of the $\myJ$ function}

Analogously to $J$ we want to estimate $\myJ$ by the ratio of
two estimators for the denominator and numerator in Definition
\ref{myJdefn}. The stationarity of $X$ and Fubini's Theorem yield
\[
\intl_W\myF(B(0,r)\cap W_{-x})\,dx
=\E\intl_W\1\{X\cap B(x, r)\cap W\neq\emptyset\}\,dx,\]
and so the denominator of (\ref{myJ}) becomes
\begin{equation}
  \label{den}
  \intl_W [1-\myF(B(0,r)\cap W_{-x})]\,dx
  =|W|\left[1-\frac1{|W|}\E|W\cap \big((X\cap W)\oplus B(0, r)\big)|\right],
\end{equation}
where $\oplus$ denotes Minkowski addition.

Applying the Campbell-Mecke formula \cite[(4.4.3)]{skm89} we find
\beqnas
\intl_W\myG(B(0,r)\cap W_{-x})\,dx
&=&\intl_W\P^{0!}(X\cap B(0, r)\cap W_{-x}\neq\emptyset)\,dx\\
&=&\frac1\lambda\E\suml_{x\in X\cap W}
         \1\{X\cap B(x, r)\setminus\{x\}\cap W\neq\emptyset\}.
\eeqnas
Let $d(x,A)$ denote the Euclidean distance from a point $x\in\R^{d}$
to a set $A\subseteq\R^{d}$. The numerator of
(\ref{myJ}) can then be expressed as
\begin{equation}
  \label{num}
  \intl_W\! [1-\myG(B(0,r)\cap W_{-x})]\,dx
  =|W|\!\left[1-\frac1{\lambda|W|}\E\!\!\!\!\suml_{x\in X\cap W}\!\!\!\!
         \1\{d(x, X\cap W\setminus\{x\})\leq r\}\right]\!.
\end{equation}
The two results (\ref{den}) and (\ref{num}) allow uncorrected
estimation of $\myJ(r)$ by 
\begin{equation}
{\myJhat}(r):=
  \frac{1-\frac1{\#(X\cap W)}\suml_{x\in X\cap W}
            \1\{d(x, X\cap W\setminus\{x\})\leq r\}}
       {1-\frac1{|W|}|W\cap \big((X\cap W)\oplus B(0,r)\big)|}
 \label{JWhat}
\end{equation}
which is the uncorrected estimate of the $J$ function referred to in the 
introduction. 

Thus the uncorrected estimate of the $J$ function, based on the
uncorrected (EDF) estimates of $F$ and $G$, can be thought of as
a ratio unbiased estimator of the $\myJ$ function. As was shown in the
previous section, the $\myJ$ function can be  interpreted in the same way as 
the $J$ function. Consequently, the uncorrected estimate of the $J$ 
function, unlike the uncorrected estimates of $F$, $G$ or $K$,  can be used 
directly as an interpretive statistic in classifying deviations from the 
Poisson process.

\section{Simulation study}

This section reports the results of simulation studies 
comparing the uncorrected estimator $\myJhat$ 
with ``corrected'' estimators $\widehat J$.
In \S\ref{S:example} we show the results for a single simulated pattern;
\S\ref{S:mean+var} reports the means and variances of the estimators of $J$
in a simulation study. These results show that $\myJhat$ typically has
smaller variance than the corrected estimators.
In \S\ref{S:test-stat} and \S\ref{S:power}
we consider the power of hypothesis tests based on the uncorrected estimator 
$\myJhat$.
It is not clear, a priori, whether ``corrected'' estimators $\widehat J$
or uncorrected estimators $\myJhat$ will yield more powerful tests.
There are two competing effects:
the variance of $\myJhat$ is smaller than that of $\widehat J$, but
$\myJ$ is less sensitive than $J$ to departures from the Poisson process
(say) according to Proposition~\ref{bound}.


For the comparisons which follow, the reduced sample (border method) 
estimator $\widehat{J}_{rs}$ was adopted.
However, for the purpose of highlighting the variations which exist between
corrected estimates of the $J$ function, the Kaplan-Meier estimator 
$\widehat{J}_{km}$ \cite[described in][]{badd:gill97} is also considered 
in many cases.

\subsection{Empirical example}
\label{S:example}
\begin{figure}[p]
\begin{center}
\includegraphics[angle=-90,scale=0.5]{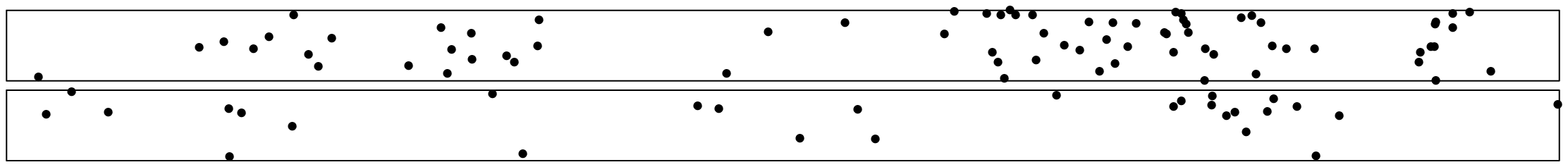}
\includegraphics[angle=-90,scale=0.25]{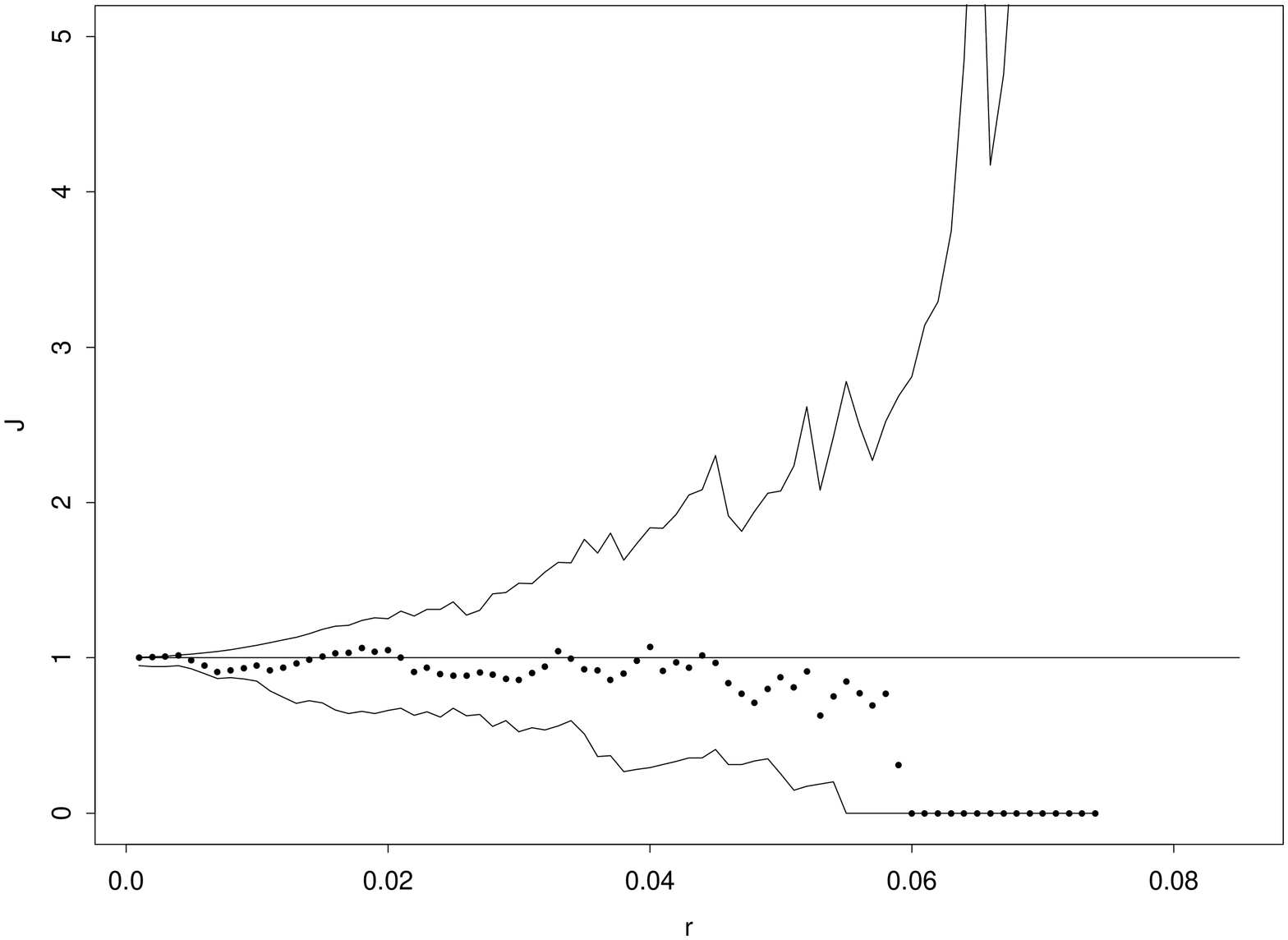}
\includegraphics[angle=-90,scale=0.25]{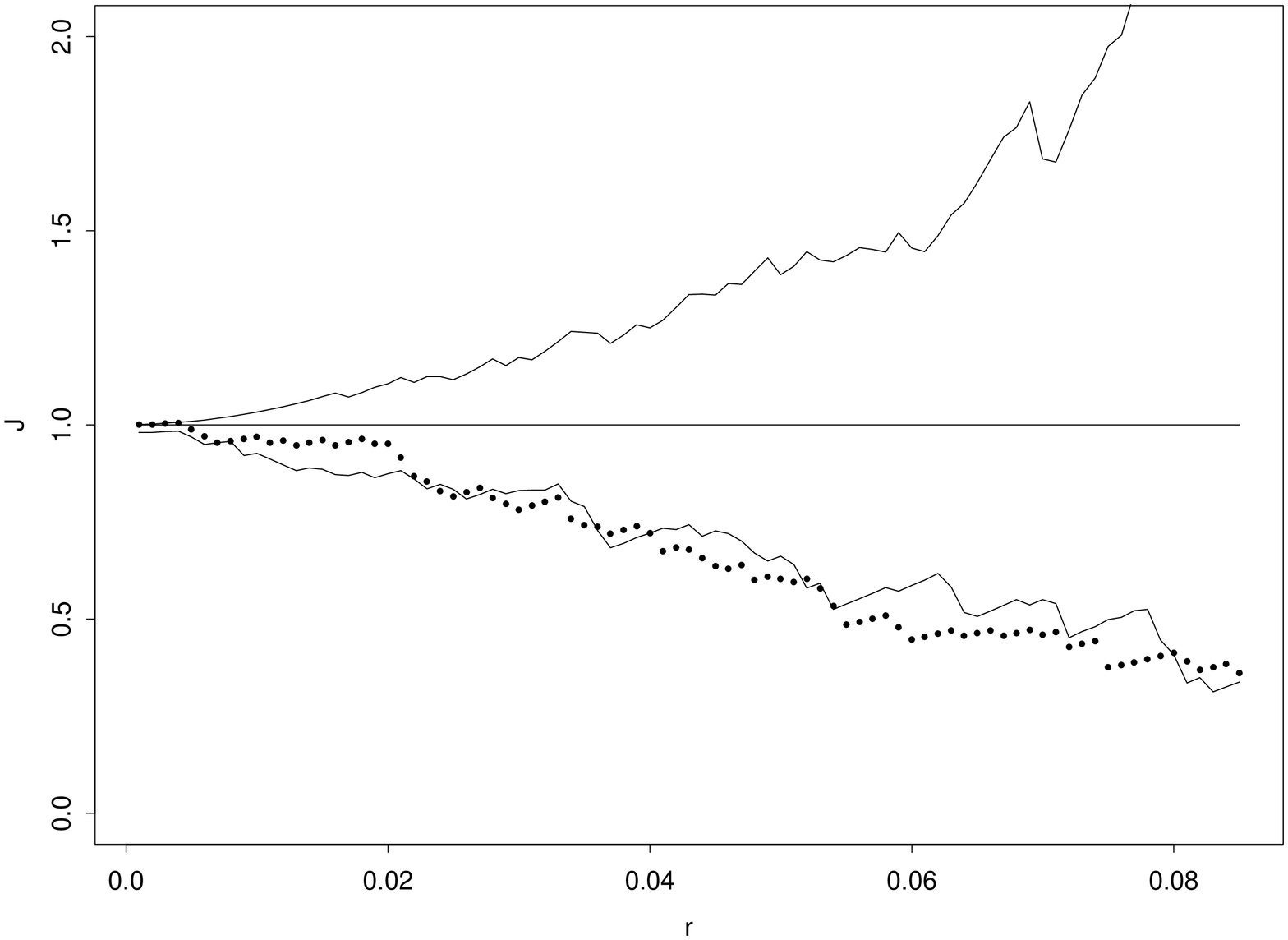}
\caption{Top: empirical data. Bottom: empirical $J_{rs}$ (left) and $J_{W}$ (right)
functions (points) and envelope of 99 simulations of a binomial process with the
same intensity (solid lines).}
\label{window}
\end{center}
\end{figure}

This example highlights the use of $\myJhat$ as a qualitative summary statistic.
Consider the point pattern given in Figure \ref{window}. This pattern is a
realisation of a Mat\'ern Cluster Process, intensity $\lambda=100$, cluster
radius $R=0.1$ and mean number of offspring $\mu=4$, observed within
an observation window consisting of two rectangular windows, 3.125 by 0.16 
units, separated by 0.02 units. (The Mat\'ern Cluster Process is a Neyman-Scott
process in which offspring are uniformly distributed in disc of radius $R$
about (Poisson) parent points. The number of offspring per parent point is
Poisson with mean $\mu$ \citep[p. 159]{skm89}).

Figure \ref{window} also displays the corresponding estimates $\myJhat$ and
$\widehat{J}_{rs}$ of the given point pattern, together with envelopes
of 99 simulations of a binomial process of the same intensity, observed within
the same window. The results in this case are marked; the uncorrected
estimate suggests strong evidence of clustering in the point pattern, while
the corrected estimate appears to suggest no evidence of clustering. Of course,
the results are not surprising given the severity of the edge effect introduced
by the ``censoring'' of the middle seventeenth of the window compared to the
relatively small bias this introduces. However, it does illustrate the
possible benefit of using $\myJhat$ in certain situations.

As the empirical use of $J_{W}$ is the same as the $J$ function, readers 
interested in further examples of the analysis of empirical data are referred 
to \cite{kersch98}, \cite{kersch:al98}, and  \cite{kersch:al97}.

\subsection{Mean and variance}
\label{S:mean+var}
\begin{figure}[p]
\begin{center}
\includegraphics[scale=0.4]{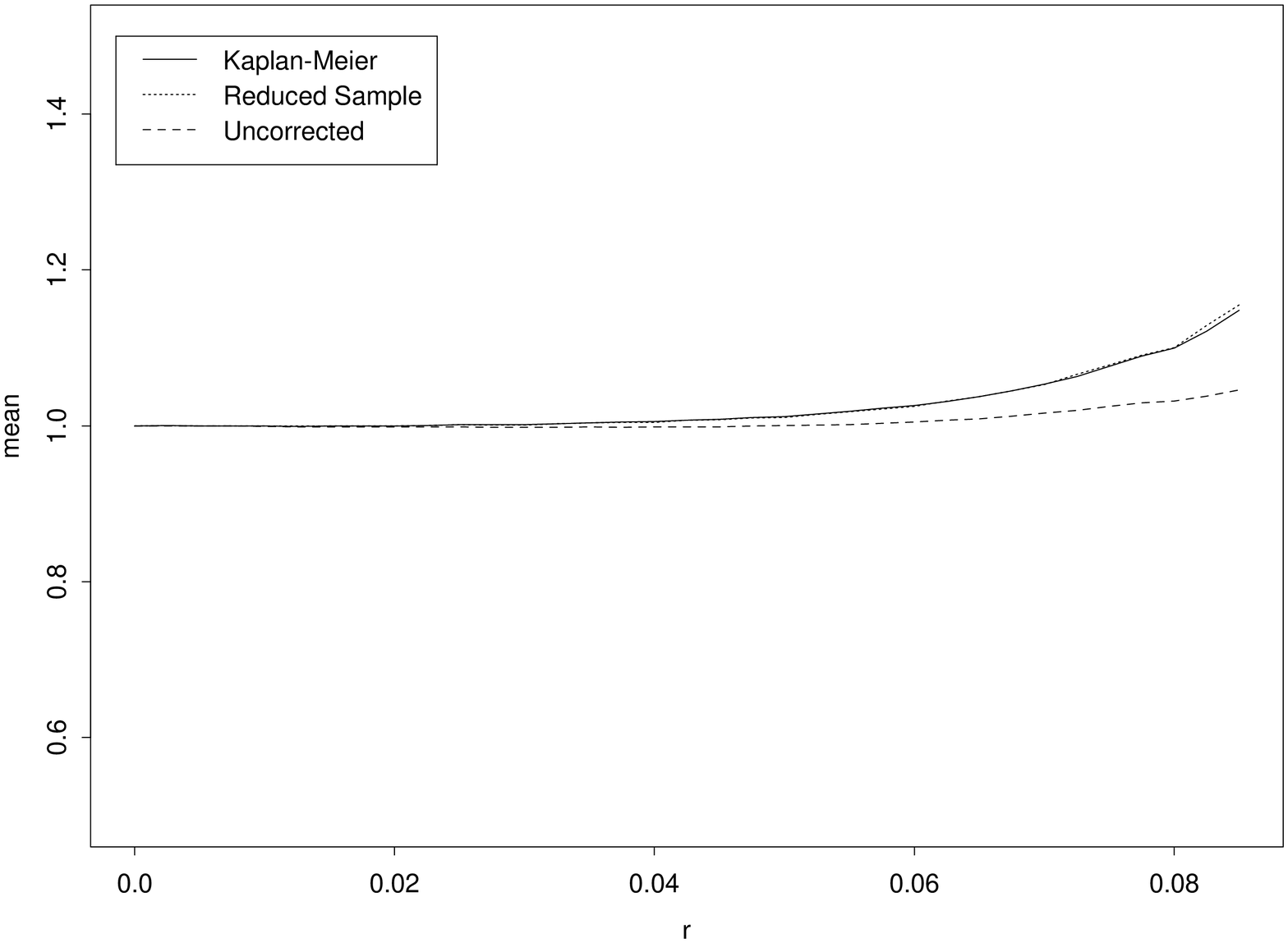}
\includegraphics[scale=0.4]{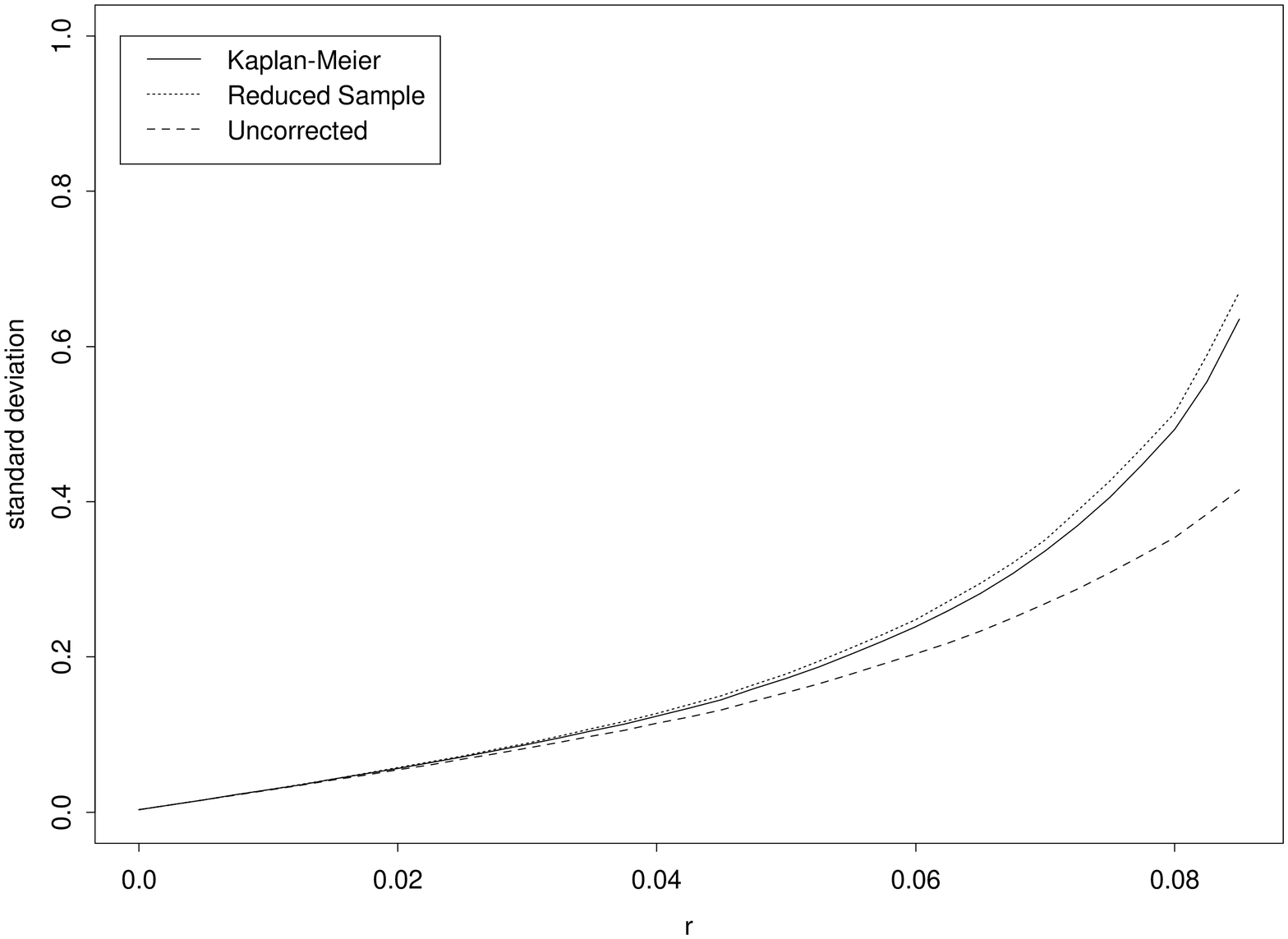}
\caption{Results in the unit square, intensity 100.
Mean (top) and standard deviation (bottom) of $J$ estimators as a 
function of $r$ for a Poisson process}
\label{meanstdev}
\end{center}
\end{figure}

>From the previous example it is clear that the uncorrected
estimate of the J function may be superior in some situations.
To examine whether this was true more generally, a number of
simulations were conducted to compare the corrected and uncorrected 
methods across a range of processes. 
This began with the estimation of the mean and standard deviation of the 
three estimators ($\myJhat$, $\widehat{J}_{km}$, and $\widehat{J}_{rs}$)
based on $10,000$ realisations of a Poisson process with intensity
$100$, in a unit square window, with the results presented in Figure
\ref{meanstdev}.

With increasing $r$, the distributions of the $\widehat{J}_.$ (that is, the estimators
$\myJhat$, $\widehat{J}_{km}$, and $\widehat{J}_{rs}$) become
skewed to the right; for large $r$ there is substantial mass above 
$\widehat{J}_.(r)=2$. As a result all three estimators are positively biased 
for large values of $r$.  Empirically it was found that a square root
transformation approximately symmetrised the distribution. As
expected, the sample standard deviation of the estimates increases with $r$,
as the denominator of each estimator decreases with $r$.
However, $\myJhat$ is less biased and has lower variance than
$\widehat{J}_{km}$ and $\widehat{J}_{rs}$. 

These simulations were repeated for two processes with more substantial
edge effects (namely, a Poisson process of intensity $\lambda=25$
in a unit window, and a Poisson process of intensity $\lambda=10$ in a
$10$ by $1$ rectangular window). In both cases the results were qualitatively
similar to those above.

\begin{figure}[p]
\begin{center}
\includegraphics[angle=-90,scale=0.25]{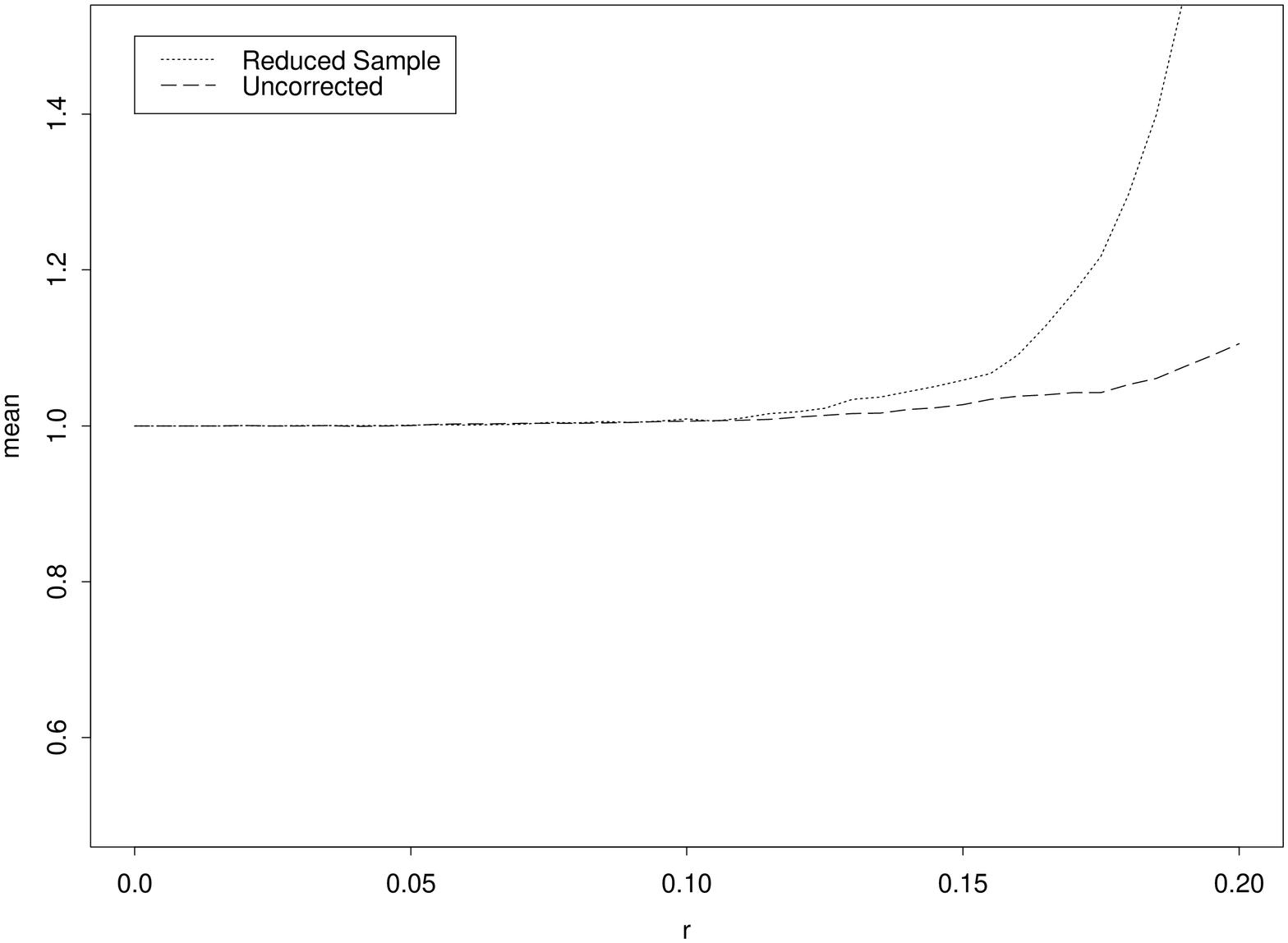}
\includegraphics[angle=-90,scale=0.25]{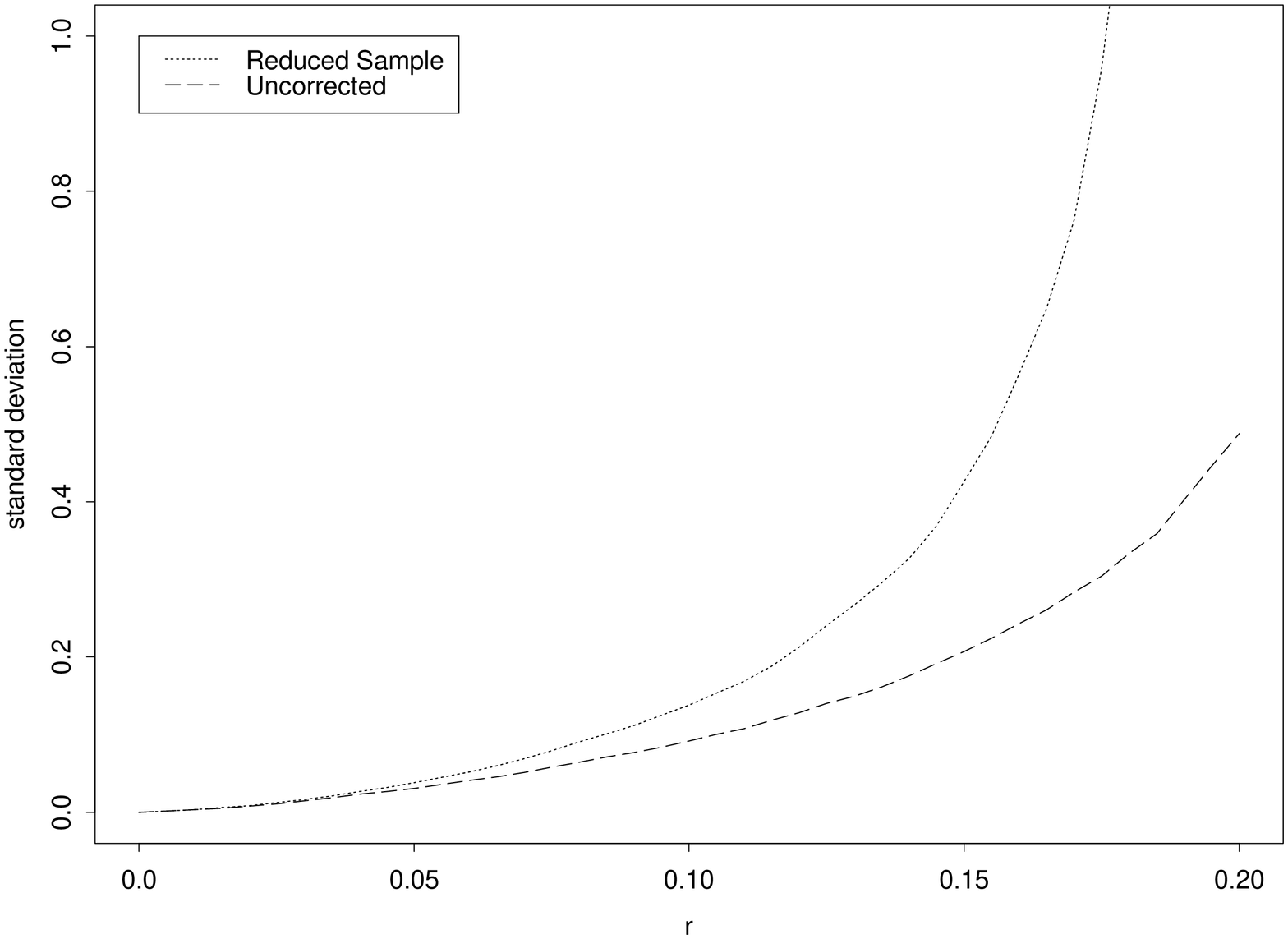}
\includegraphics[angle=-90,scale=0.25]{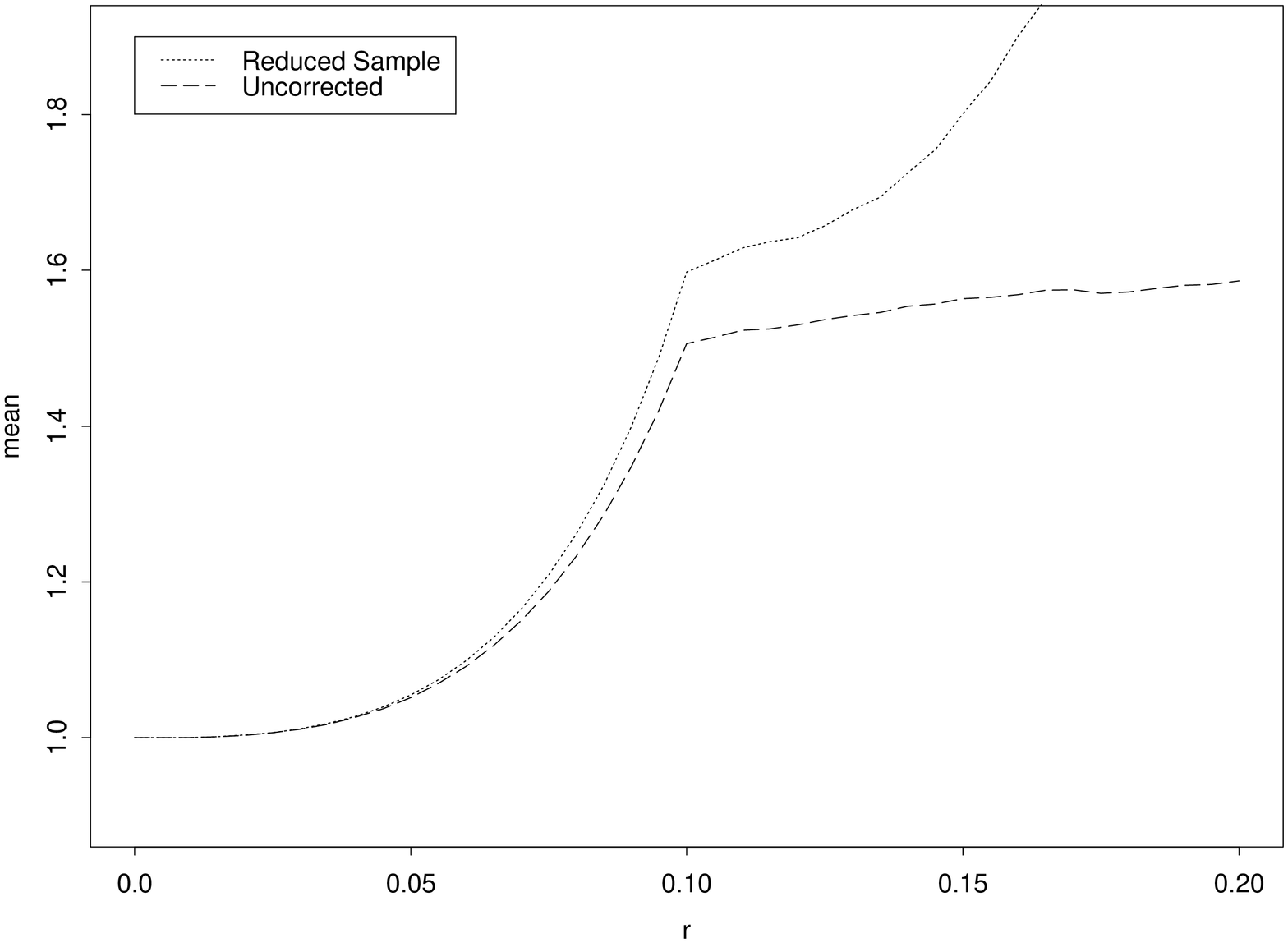}
\includegraphics[angle=-90,scale=0.25]{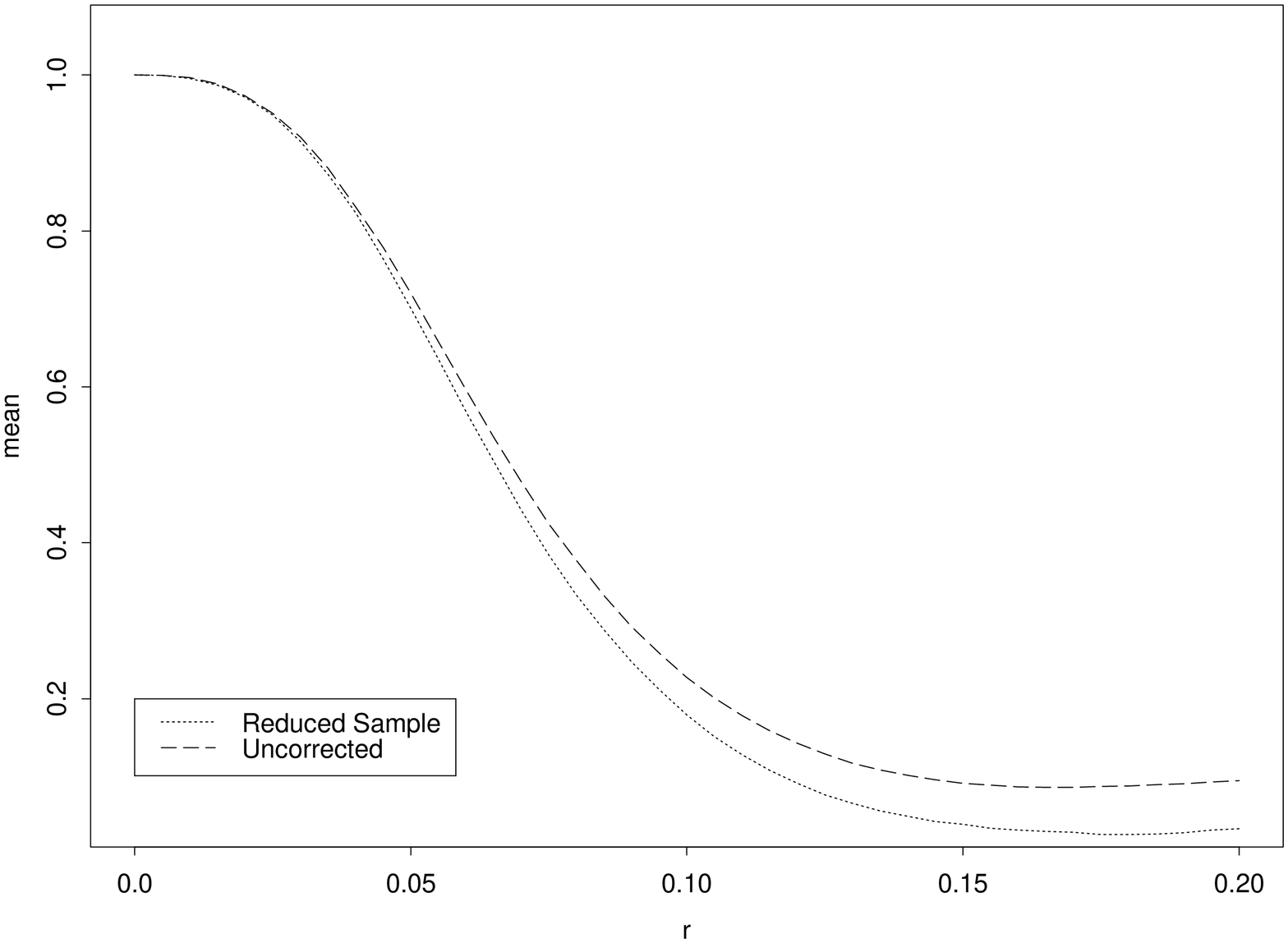}
\caption{Results in the unit cube, intensity=100.
Top: Mean (left) and standard deviation (right) of $J$ estimators as a 
function of $r$ for a Poisson process.
Bottom left: Mean of $J$ estimators for a Mat\'ern model II process ($R=0.1$).
Bottom right: Mean of $J$ estimators for a Mat\'ern cluster process ($\mu=4$,
$R=0.1$). }
\label{cube}
\end{center}
\end{figure}

In addition, further simulations were conducted for point patterns in $\R^3$. 
Estimates of the means and standard 
deviations of $\widehat{J}_{km}$ and $\widehat{J}_{rs}$
based on $1000$ realisations in a unit cube 
were compared for the Poisson process and two alternatives:
Mat\'ern 
hard-core \citep[p. 163]{skm89}  and  Mat\'ern cluster processes 
for a range of parameter values. Some of the key results are presented in 
Figure \ref{cube}.
As expected, $\myJhat$ is reliable over a wider domain than
$\widehat{J}_{rs}$. For hard-core processes, the standard deviation of
$\widehat{J}_{rs}$ is considerably bigger than that of $\myJhat$ and the 
difference grows with the hard-core radius. For cluster processes the 
differences are far less apparent, however the overall tendency of
$\myJhat$ to have lower variance is also confirmed for
this class of processes. Note also that with both alternative processes
the mean of $\myJhat$ is bounded by $1$ and $\widehat{J}_{rs}$, as expected.

It is interesting to note that, unlike the $J$ function estimators, the domain
of the $\myJ$ function estimator for a given point process realisation can be
easily calculated. The $\myJ$ estimator is defined for all
$r<r_{F_{max}}$, where $r_{F_{max}}$ is the maximum nearest-point
distance (the largest distance, over all points in the window, from a
point to the nearest point of the process). Also $\myJ(r)=0$ for any
$r_{G_{max}}\leq r<r_{F_{max}}$ if $r_{G_{max}}<r_{F_{max}}$, where
$r_{G_{max}}$ is the maximum nearest-neighbour distance (the largest
distance, over all points of the process within the window, from a point of the 
process to the nearest other point of the process).  The value $r_{F_{max}}$
is however an upper bound on the domain of both the Reduced Sample and
Kaplan-Meier estimators.

\subsection{The test statistic}
\label{S:test-stat}
We now aim to compare the power of the $J_W$ function with the edge
corrected estimators of the $J$ function in testing the Poisson
hypothesis in the two dimensional case. We restricted ourselves to
this estimation problem in view of the problems with estimating the
range of interaction using the $J$ function reported in
\cite{kersch:al98}. 

The distribution of the following test statistic for each of the three
estimators was estimated:
\begin{equation}
\tau=\int_{0}^{r_{0}}\frac{\widehat{J}_.(r)-1}{\widehat{\sigma}(r)}\,dr,
\end{equation}
where $\widehat{\sigma}$ denotes the sample standard deviation of
$\widehat{J}_.(r)$ under the Poisson hypothesis.  This form of
test statistic was chosen, as opposed to a squared integrand, because
of the skewed nature of the distributions of $\widehat{J}_.$.  The
distributions of the test statistics were estimated by 
a discrete sum and based on $10,000$ realisations of a Poisson
process.  The upper limit of integration $r_{0}$ was chosen to be the
0.9 quantile of the $F$ function (for intensity $100$, $r_{0}\approx
0.856$).  Having estimated the distribution, the $0.025$ and $0.975$
quantiles were obtained for use in a two-sided $5\%$ significance test
for deviation from a Poisson process.  One-sided $5\%$ significance
tests were also constructed to test for clustering or regularity by
considering the $0.05$ and $0.95$ quantiles, respectively.


\subsection{Power of tests using the various estimators}
\label{S:power}

\begin{figure}[p]
\begin{center}
\includegraphics[scale=0.4]{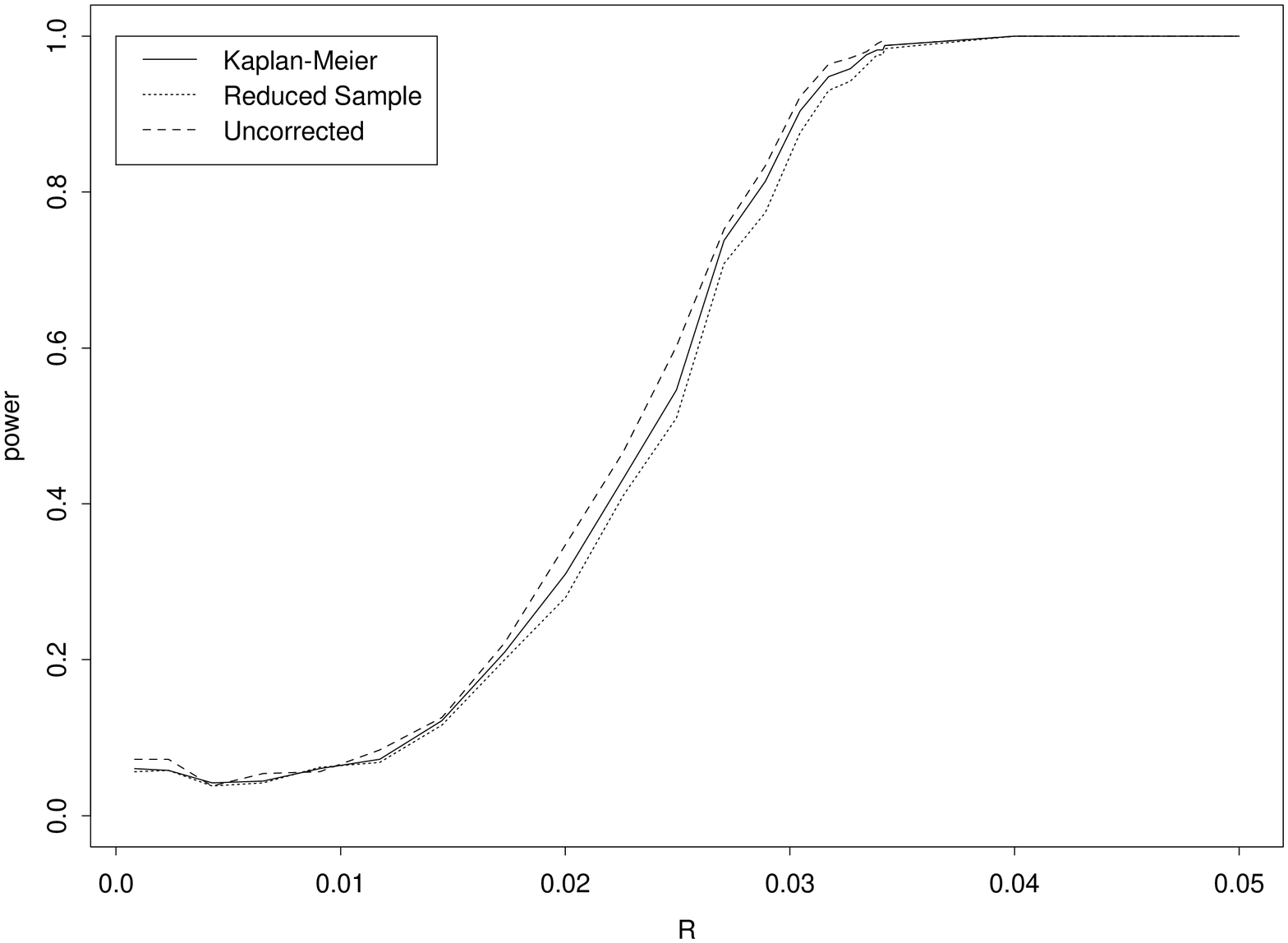}
\includegraphics[scale=0.4]{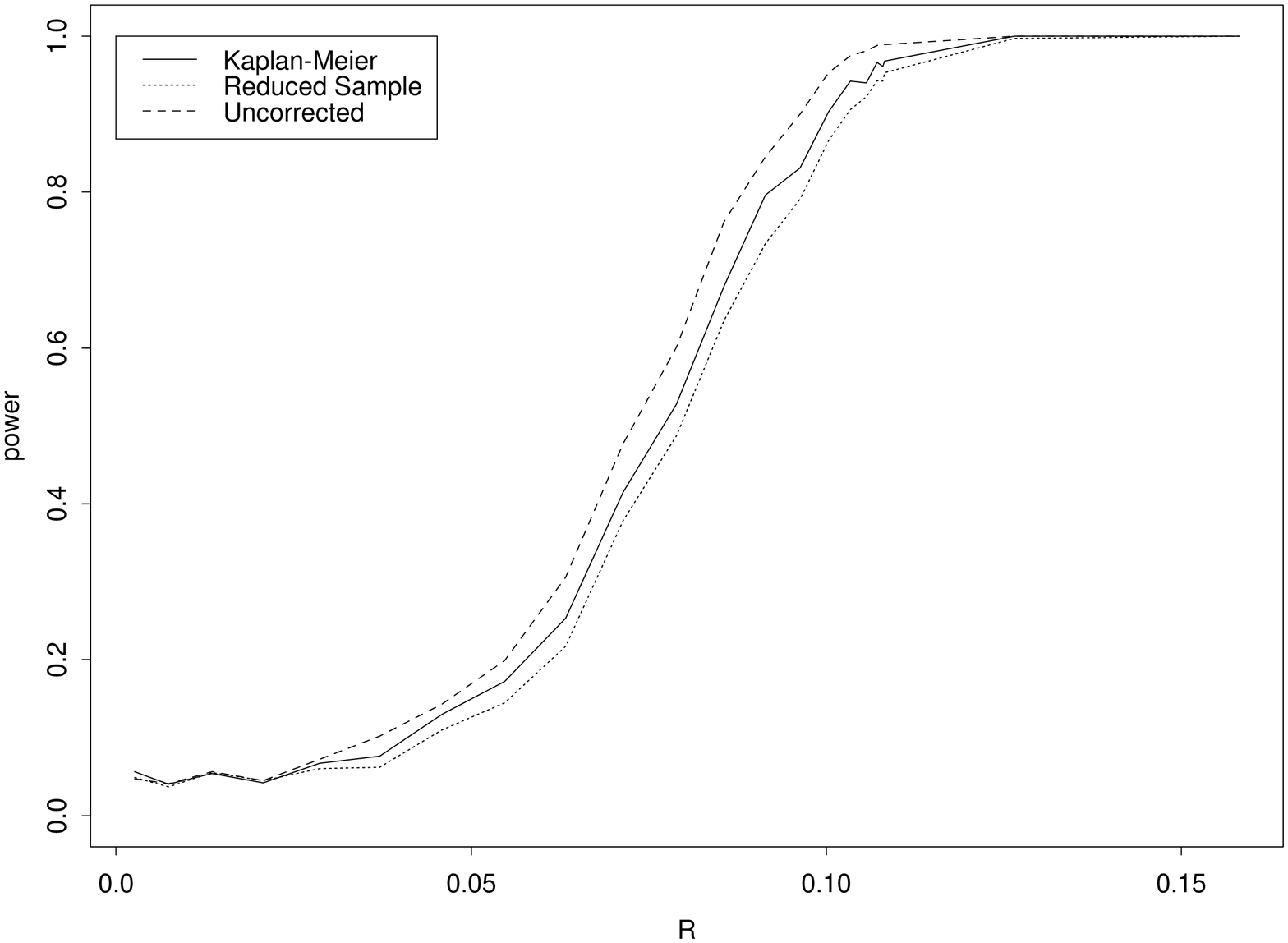}
\caption{Power against Mat\'ern Model II as a function of  hard-core 
radius $R$.
Top: unit square, intensity 100. 
Bottom: $[0,1]\times[0,10]$ rectangle, intensity 10}
\label{MMII}
\end{center}
\end{figure}

In order to estimate the power of the hypothesis test described above, 
realisations from alternative point processes were generated and the proportion 
of the realisations rejected by the hypothesis test was recorded. The first
class of point processes considered was the Mat\'ern hard-core
process, with hard-core radius $R$. For each of $22$ values of $R$,
$1000$ realisations were generated. The proportion of rejections is
presented in Figure \ref{MMII}. Note that as $R\rightarrow 0$ the
model approaches the Poisson process, so we expect all power curves to
approach $0.05$ $(5\%)$ as $R\rightarrow 0$.  All three
estimators have very similar power curves, with the $\myJ$ estimator
at least as powerful as the two $J$ function estimators for all values
of $R$.

\begin{figure}[p]
\begin{center}
\includegraphics[scale=0.4]{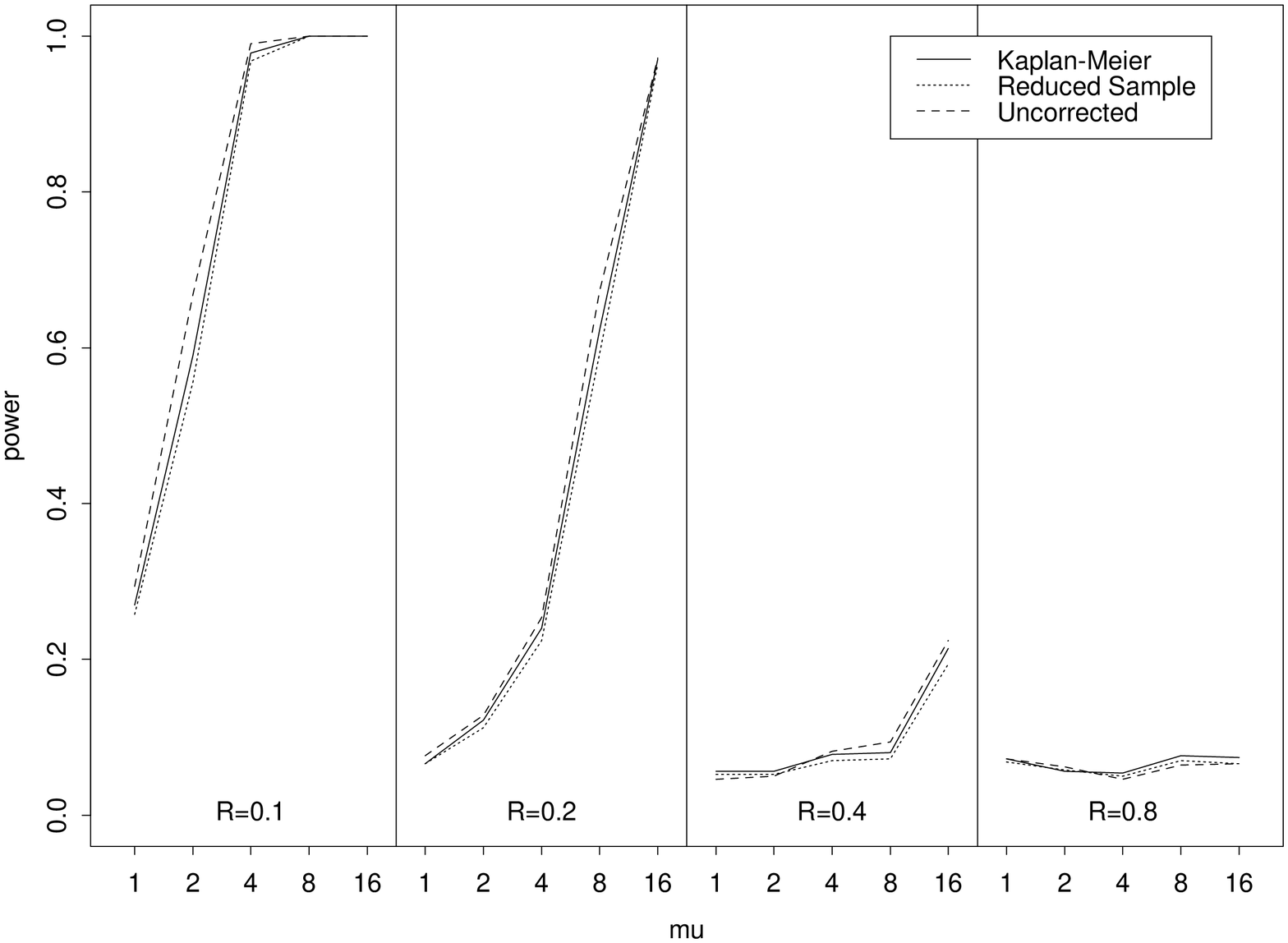}
\includegraphics[scale=0.4]{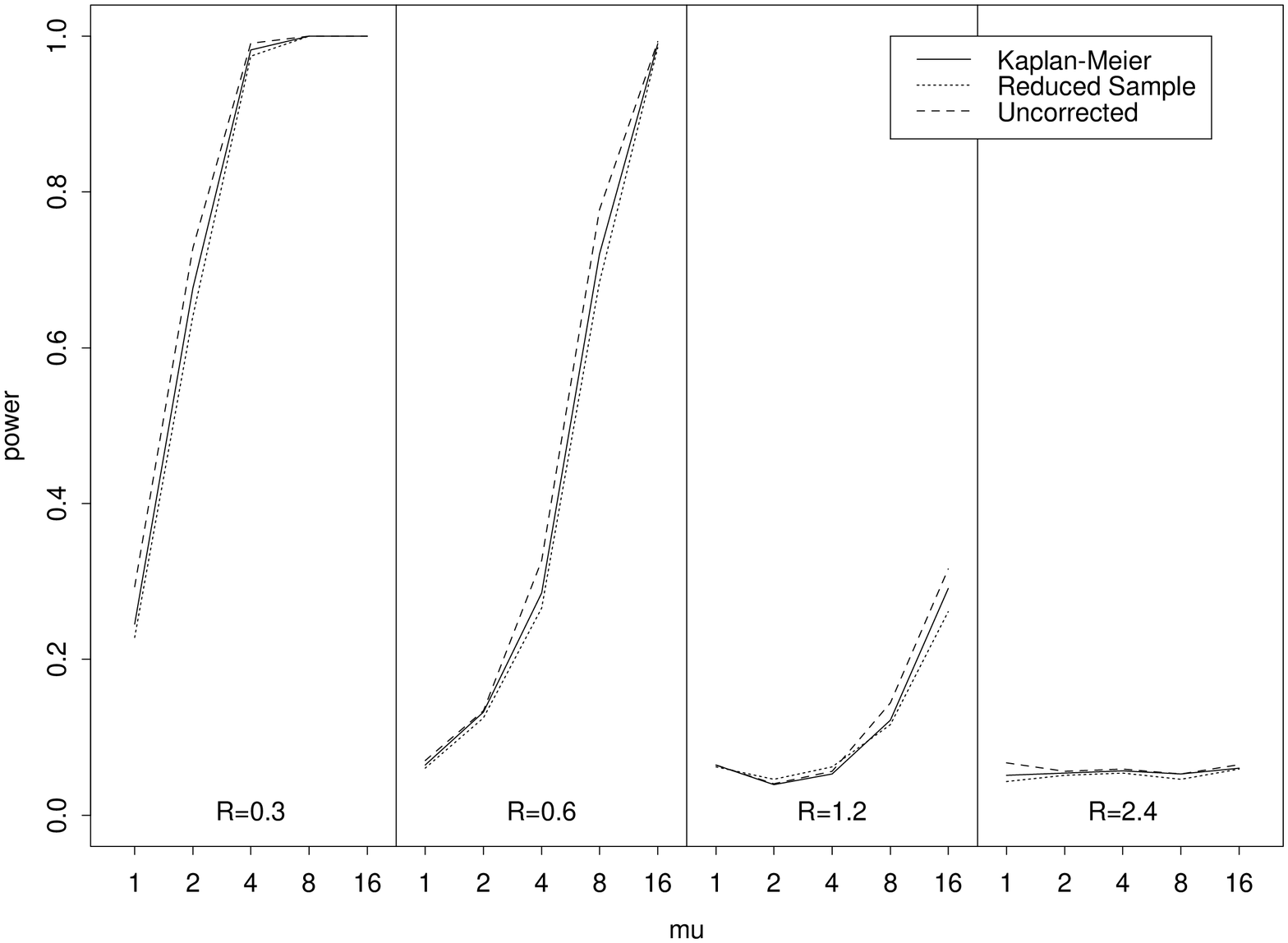}
\caption{Power against Mat\'ern Cluster Process as a function of mean cluster 
size for various cluster radii $R$.
Top: unit square, intensity 100. 
Bottom: $[0,1]\times[0,10]$ rectangle, intensity 10}
\label{MNS}
\end{center}
\end{figure}

The other class of alternative point processes considered was
Mat\'ern's cluster process. A grid of $(R,\mu)$ values was constructed
and $1000$ realisations were obtained from the corresponding Mat\'ern
cluster process. The proportion of rejections for each $(R,\mu)$ value
is presented in Figure \ref{MNS}. Once again, the curves are very
similar, with all three tests performing almost identically. Similar
results were obtained for the respective one-sided tests in both
cases as well as for the lower intensity $25$.  

The power tests for the $10$ by $1$ window support the argument that
edge effects are stronger when the boundary is relatively longer. The
resulting power function estimates against the Mat\'ern Model II and
Mat\'ern cluster process models are also presented in Figures
\ref{MMII} and \ref{MNS}, respectively. 

One important observation made while conducting these numerical
simulations was that the choice of test statistic had far more impact
on the power of the resulting hypothesis test than the choice of $J$
function estimator. For a comparison of various test statistics of
the $J$ function see \cite{thon:lies99}.

\newcommand{\noopsort}[1]{} \newcommand{\printfirst}[2]{#1}
  \newcommand{\singleletter}[1]{#1} \newcommand{\switchargs}[2]{#2#1}

\end{document}